\documentclass[12pt,twoside,a4paper]{article}
\usepackage[utf8]{inputenc}
\usepackage{mathtools}
\usepackage{physics}
\usepackage[none]{hyphenat}

\usepackage{graphicx} 
\usepackage{etoolbox}

\usepackage{amssymb,amsmath,amsthm,amsfonts,latexsym} 
\usepackage{correctmathalign}
\usepackage[mathscr]{euscript}

\usepackage{mathtools}
\usepackage[none]{hyphenat}
\DeclarePairedDelimiter{\ceil}{\lceil}{\rceil}
\DeclarePairedDelimiter\floor{\lfloor}{\rfloor}

\makeatletter
\DeclareRobustCommand*{\bfseries}{%
  \not@math@alphabet\bfseries\mathbf
  \fontseries\bfdefault\selectfont
  \boldmath
}
\makeatother

\usepackage[labelsep=quad,indention=10pt,format=hang]{caption}
\usepackage{subcaption}

\usepackage[pdftex,bookmarks,colorlinks=false]{hyperref}

\usepackage{geometry}
\usepackage[auth-lg,affil-sl]{authblk}
\usepackage{float} 
\usepackage{dpfloat}
\usepackage{multirow} 
\usepackage{hhline} %
\usepackage{booktabs} 
\usepackage{longtable} 
\usepackage{diagbox} 

\usepackage{cite}
\usepackage{url}
\usepackage{breakcites} 
\usepackage{breqn} 
\allowdisplaybreaks

\usepackage{cite}

\usepackage{color} 
\usepackage{colordvi}
\usepackage{colortab}
\usepackage{colortbl}
\usepackage{colorweb}

\usepackage{bullcntr}%
\usepackage{doi}
\usepackage{epic}
\usepackage{epigraph}

\usepackage{enumitem}

\usepackage{calc,tikz}
\usetikzlibrary{}
\usetikzlibrary{intersections, arrows, arrows.meta, automata, er, calc, mindmap, folding, positioning,patterns, decorations.markings, fit, 
shapes, matrix, positioning, shapes.geometric, through}
\tikzstyle{vertex}=[circle, draw, inner sep=1pt, minimum size=4pt]
\newcommand{\vertex}{\node[vertex]}

\tikzstyle{ann} = [fill=white,font=\footnotesize,inner sep=1pt]
\tikzstyle{arrow} = [thick,<-->,>=stealth]

\newcommand{\noi}{\noindent}

\newcommand{\Mod}[1]{\ (\mathrm{mod}\ #1)}
\newtheorem{theorem}{Theorem}[section]
\newtheorem{definition}[theorem]{Definition}

\newcommand{\keywordsname}{Keywords}
\newenvironment{keywords}{
	\normalfont
	\list{}{
		\labelwidth0pt
		\leftmargin0pt \rightmargin\leftmargin
		\listparindent\parindent \itemindent0pt
		\parsep0pt
		}
	\item[\hskip\labelsep\bfseries\keywordsname:] \small}{\endlist}

\newcommand{\mscname}{MSC 2020}
\newenvironment{msc}{
	\normalfont
	\list{}{
		\labelwidth0pt
		\leftmargin0pt \rightmargin\leftmargin
		\listparindent\parindent \itemindent0pt
		\parsep0pt
		}
	\item[\hskip\labelsep\bfseries\mscname:] \small}{\endlist}

\definecolor{ududff}{rgb}{0.30196078431372547,0.30196078431372547,1}
\definecolor{xdxdff}{rgb}{0.49019607843137253,0.49019607843137253,1}

\newtheoremstyle{casesty}         
{}                   
{}                   
{\upshape}           
{}                   
{\itshape}         
{:-}                   
{1em}                
{}                   
\theoremstyle{casesty}
\newtheorem{case}{Case}

\newtheoremstyle{subcasesty}         
{}                   
{}                   
{\upshape}           
{}                   
{\itshape}          
{:-}                   
{1em}                
{}                   
\theoremstyle{subcasesty}

\newtheoremstyle{schemesty}         
{}                   
{}                   
{\upshape}           
{}                   
{\itshape}          
{:-}                   
{1em}                
{}                   
\theoremstyle{schemesty}

\title{\sc Equitable Dominator Coloring of Line Graphs of Some Graphs}

\author{{\bf Phebe Sarah George$^\ast$ and Sudev Naduvath$^\dag$}}
\affil{\tt Department of Mathematics\\\tt  CHRIST (Deemed to be University) \\ \tt Bangalore-560029, India.\\
$^\ast${\tt phebe.george@res.christuniversity.in} \\
$^\dag${\tt sudev.nk@christuniversity.in}}
\date{}

\begin{document}
\maketitle
\hrule

\begin{abstract}

A proper vertex coloring of a graph $G$ such that every vertex of $G$ dominates at least one color class and the cardinalities of the color classes differ by at most $1$ is called an equitable dominator coloring of $G$. The minimum number of colors used in an equitable dominator coloring of a graph $G$ is called the equitable dominator chromatic number of $G$ and is denoted by $\chi_{ed}(G)$. This article explores the concept of equitable dominator coloring for the line graph $L(G)$ of some elementary graph classes. 

\vspace{0.2cm}

\vspace{0.4cm}

\hrule

\end{abstract}

\section{Introduction}
For all the terms and definitions in graph theory, one may refer to \cite{harary1,dbw1}. For further topics in graph coloring and domination in graphs, refer to \cite{cz1,fundamentalsofdom, topicsindom,mk1}.  Unless mentioned otherwise, all graphs under consideration are undirected, simple, connected, and finite.

The \emph{line graph} of a graph $G$, denoted by $L(G)$, is a graph with $V(L(G))=E(G)$ such that two vertices $e_i,e_j$ of $L(G)$ are adjacent if the corresponding edges have a vertex common in $G$. Note that we use the same label of an edge in $G$ to denote the corresponding vertex in $L(G)$.

 A \emph{graph coloring} of a graph $G$ is an assignment of colors, weights, or labels to the elements - vertices, edges, or faces - of $G$. In this paper, the term graph coloring refers to the assignment of colors to the vertices of $G$.  A graph $G$ is said to be \emph{properly colored} when no two adjacent vertices of $G$ are assigned the same color. The minimum number of colors required to properly color a graph $G$ is called the chromatic number of $G$, and is denoted by $\chi(G)$. The colorings considered in this paper are all proper vertex colorings of a graph $G$. A \emph{color class} $V_{j}\subseteq V(G)$ of a graph $G$ is the set of all vertices of $G$ that receive the same color $c_j$ in a coloring $c$ of $G$. 

 A set $S \subseteq V(G)$ is said to be a \emph{dominating set} of a graph $G$ if for all $v \in V(G)$  either $v \in S$ or $v\sim u$ such that $u\in S$.  The \emph{domination number} of a graph $G$, denoted by $\gamma(G)$, is the cardinality of a minimum dominating set of a graph $G$.

Combining the concepts of domination and coloring, the concept of dominator coloring of a graph $G$ was introduced in \cite{DRsafeclique} as follows: 
\begin{definition}\label{Domcol}{\rm \cite{DRsafeclique}
A \emph{dominator coloring} of a graph $G$ is a coloring of $G$ such that every vertex in the vertex set of $G$ dominates all the vertices of at least one color class, possibly its own color class. The \emph{dominator chromatic number} of $G$, denoted by $\chi_{d}(G)$, is the minimum number of colors used in a dominator coloring of a graph $G$.
}\end{definition}
A variant of proper vertex coloring of a graph $G$ called the equitable coloring of graphs was introduced in  \cite{1973equitable} as follows.

\begin{definition}\label{equi}{\rm \cite{1973equitable}
An equitable coloring of a graph $G$ is a proper vertex coloring of $G$ such that the cardinalities of any two color classes differ by at most 1. The minimum number of colors used in an equitable coloring of a graph $G$ is called the \emph{equitable chromatic number} of $G$, denoted by $\chi_{e}(G)$.
}\end{definition}

\noi A new variant of domination-related coloring, called the equitable dominator coloring of a graph $G$, was introduced in \cite{eqdom1} as follows.

\begin{definition}{\rm \cite{eqdom1}
An \textit{equitable dominator coloring} of a graph $G$ is a proper coloring of $G$ such that every vertex in $V(G)$ is assigned a color such that it dominates at least one color class, possibly its own color class, and the cardinalities of the color classes differ by at most one.
The minimum number of colors used in an equitable dominator coloring of $G$ is called the \emph{equitable dominator chromatic number} of $G$, and it is denoted by $\chi_{ed}(G)$.}
\end{definition}

The concept of equitable dominator coloring of graphs was introduced in \cite{eqdom1} and studied for several families of graphs, including path graphs, cycle graphs, and cycle-related graphs such as wheel graphs, helm graphs, and closed helm graphs. The equitable dominator chromatic number for the complements of these graph classes was also obtained. Motivated by the above-mentioned studies, in this paper, we obtain an equitable dominator coloring and the corresponding equitable dominator chromatic number for the line graph of the above-mentioned graph families.

\section{Equitable Dominator Chromatic Number of Line Graphs of Some Graphs}

 Note that  $L(P_n) \cong P_{n-1}$, and $L(C_n) \cong C_n$. Since the equitable dominator coloring and the corresponding equitable dominator chromatic number of these graph classes have been obtained in \cite{eqdom1}, we initiate our present study by exploring the concept of equitable dominator coloring in the line graph of a bi-star graph.

A \emph{bi-star} denoted by $S_{a,b}$; $a,b\geq 2$, is a graph obtained by joining the two universal vertices, say $u$ and $v$, of two star graphs $K_{1,a}$ and $K_{1,b}$, respectively, by an edge. 

\begin{theorem}
For $a,b\geq 2$, $\chi_{ed}(L(S_{a,b}))=\max\{a,b\}+1$    
\end{theorem}

\begin{proof}
Let $u_i;1 \leq i \leq a$ be the pendant vertices of the star graph $K_{1,a}$, and $v_i;1 \leq i \leq b$ be the pendant vertices of the star graph $K_{1,b}$. Let $u,v$ be the universal vertices of the two star graphs $K_{1,a}$ and $K_{1,b}$ respectively. Let $e_i;1 \leq i \leq a$ be the edges between the vertices $uu_i;1 \le i \le a$, and $e_i'; 1 \le i \le b$ be the edges between the vertices $vv_i;1\le i\le b$ in $S_{a,b}$. Let $e$ represent the edge joining the universal vertices $u,v$ of the two star graphs $K_{1,a}$ and $K_{1,b}$, respectively in $S_{a,b}$.
In $L(S_{a,b})$, the edge $e$ joining the universal vertices $u,v$ of the two star graphs $K_{1,a}$ and $K_{1,b}$ in $S_{a,b}$ becomes an universal vertex; hence the vertex $e$ in $V(L(S_{a,b}))$ is assigned a unique color. The vertices $\{e_i: 1\le i \le a\}$, and the vertices $\{e_i': 1\le i \le b\}$ in $L(S_{a,b})$, each induce a clique of order $a,b$ respectively, in $L(S_{a,b})$. Without loss of generality, let $\min\{a,b\}=a$. Since the edges $e_i$ and $e_i'$ are not adjacent to each other, $\omega(L(S_{a,b}))=b+1$. Thus, we have $\chi_{ed}(L(S_{a,b})) \geq b+1$. Consider a coloring $c$ of $L(S_{a,b})$ such that $c(e_i')=c_i; 1\leq i \leq b$, and $c(e_i)=c(e_i');1 \le i \le a$. Let the universal vertex $e$ in $L(S_{a,b})$ be assigned the color $c_{b+1}$. Here, the cardinality of every color class in the coloring $c$ of $L(S_{a,b})$ is at most 2, and all the vertices of $L(S_{a,b})$ dominate the color class $V_{b+1}$.  Thus, $c$ is an equitable dominator coloring of $L(S_{a,b})$ and the result follows. 
\end{proof}

\begin{theorem}
For $a,b\geq 1$ , $\chi_{ed}(L(K_{a,b}))=\max\{a,b\}$.

\end{theorem}

\begin{proof}
 The line graph of a complete bipartite graph $K_{a,b}$ is a graph of order $ab$ such that $(i,j)$ denotes a vertex of $L(K_{a,b})$, where $1 \leq i \leq a,1 \leq j \leq b$ and every vertex $(i,j)$ represents an edge connecting the vertices $a_i$ and $b_j$ in the original graph $K_{a,b}$. In $L(K_{a,b})$, two vertices $(i,j)\sim (i',j')$ if and only if $i=i'$ or $j=j'$ but not both.

\textit{Claim 1:} $\omega(L(K_{a,b}))=\max\{a,b\}$.
From the adjacency of vertices $(i,j)$ in $L(K_{a,b})$, we can see that for a fixed $i$ and $1 \leq j \leq b$. the vertices $(i,j)$ form a clique of size $b$. Similarly, for a fixed $j$ and $1 \leq i \leq a$, the vertices $(i,j)$ form a clique of size $a$. If exists, take a clique of size greater than $\max\{a,b\}$. Without loss of generality, let $a < b$. If $\omega(L(K_{a,b}))=b+1$, this implies for a fixed $i$ and $1 \leq j \leq b$ the set of vertices $(i,j)$ which forms a clique of size $b$ is adjacent to at least one vertex $(i',j')$ where $1 \leq i' \leq a$ but $i' \neq i$ and $1 \leq j' \leq b$.

\noi Here, we have the following cases.

\begin{case}
For $i \neq i'$, consider $j'=j$ for a single value of $j$. Then, the vertex $(i,j)\sim (i',j')$ only for the single value of $j=j'$. The vertex $(i',j')$ will not be adjacent to any other vertex $(i,j)$, for all the other values of $j$. Hence, a clique of order $b+1$ is not possible.
\end{case}

The same argument holds for all cliques of order greater than or equal to $b+1$. Hence, $\omega(L(K_{a,b}))=\max\{a,b\}$. 

\textit{Claim 2:} $\chi(L(K_{a,b}))=\max\{a,b\}$. We have a result in \cite{edgechrovertex} that for a graph $G$ the edge chromatic number $\chi'(G)$ is equal to the vertex chromatic number of the line graph of the graph $G$. It is also proved that the edge chromatic number of $K_{a,b}$ is $\max\{a,b\}$. Thus, combining both the results, $\chi(L(K_{a,b}))=\max\{a,b\}$.

\textit{Claim 3:} $\chi_{ed}(L(K_{a,b}))=\max\{a,b\}$. 
Consider a coloring $c$ of $L(K_{a,b})$ such that $c(i,j)=c_{i+j-1}$; for $1 \leq i \leq a, 1\leq j \leq b$ and the suffixes are taken under addition modulo $b$. By the adjacencies defined in \cite{hoffman1964line},  and the coloring mentioned we see that $\chi_{ed}(L(K_{a,b}))\leq \max\{a,b\}$. Since, $\omega(L(K_{a,b}))=\max\{a,b\}$, $\chi_{ed}(L(K_{a,b}))\geq \max\{a,b\}$.  Hence the result.
\end{proof}

\noi Recall that a \emph{wheel graph} denoted by $W_{1,t}; t 
\geq 3$ is a graph obtained by the join of $K_1+C_t$.

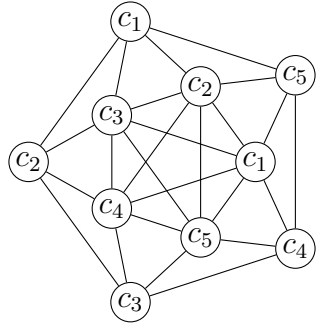
\begin{figure}[h]
\centering
\begin{tikzpicture}[scale=0.3] 
\vertex (1) at (0*360/5:3.5) []{$c_1$};
\vertex (2) at (1*360/5:3.5) []{$c_2$};
\vertex (3) at (2*360/5:3.5) []{$c_3$};
\vertex (4) at (3*360/5:3.5) []{$c_4$};
\vertex (5) at (4*360/5:3.5) []{$c_5$};
\vertex (8) at (1*360/10:6.5) []{$c_5$};
\vertex (10) at (3*360/10:6.5) []{$c_{1}$};
\vertex (12) at (5*360/10:6.5) []{$c_{2}$};
\vertex (14) at (7*360/10:6.5) []{$c_{3}$};
\vertex (16) at (9*360/10:6.5) []{$c_{4}$};
\path
(1) edge (2)
(2) edge (3)
(3) edge (4)
(4) edge (5)
(5) edge (1)
(1) edge (3)
(1) edge (4)
(2) edge (4)
(2) edge (5)
(3) edge (5)
(1) edge (16)
(16) edge (5)
(5) edge (14)
(14) edge (4)
(12) edge (3)
(12) edge (4)
(3) edge (10)
(10) edge (2)
(2) edge (8)
(8) edge (1)
(8) edge (10)
(12) edge (10)
(12) edge (14)
(14) edge (16)
(16) edge (8)
;
\end{tikzpicture}
\caption{Equitable dominator coloring of $L(W_{1,5})$.}
\label{wheel}
\end{figure}

\begin{theorem}
For $t \geq 3$, $\chi_{ed}(L(W_{1,t}))=t.$
\end{theorem}
\begin{proof}
Let $V(W_{1,t})=\{v\}\cup\{v_i: 1\le i\le t\}$ such that $\deg(v_i)=3; 1\le i\le t$ and $\deg(v)=t$. Let $e_i=vv_i; 1 \le i \le t$ be the spokes and $e_i'=v_iv_{i+1};1\le i\le t$ be the rim edges of the wheel graph $W_{1,t}$; where the suffixes are taken under addition modulo $t$. Since all the spokes are incident to the universal vertex $v$ in $W_{1,t}$, the vertices $e_i;1 \leq i \leq t$ of $L(W_{1,t})$ induce a clique of order $t$ in $L(W_{1,t})$. Thus, $\chi_{ed}(L(W_{1,t}))\geq t$. Let $c:V(L(W_{1,t})\xrightarrow{}\mathcal{C}$ be a coloring such that for $1\le i \le t$, $c(e_i)=c_i$ and $c(e_i')=c_{i-1}$. Note that the suffixes are taken under addition modulo $t$. Here, every vertex $e_i$ dominates the color class $\{e_{i-1},e_i'\}$ and every vertex $e_i'$ dominates the color class $\{e_i,e_{i+1}'\}$. Also note that the cardinality of every color class $V_i; 1\le i\le t$, is 2. Thus, the coloring $c$ is an equitable dominator coloring of $L(W_{1,t})$ using $t$ colors and it follows that $\chi_{ed}(L(W_{1,t}))=t$. For illustration, refer to Figure \ref{wheel}.
\end{proof}

A \textit{helm graph} denoted by $H_{1,t,t}; t \geq 3$ is a graph obtained by adjoining a pendant edge to each vertex of degree $3$ in a wheel graph $W_{1,t}$. The vertex of degree $t$ in $H_{1,t,t}$ is called the central vertex.

\begin{theorem}
For $t \geq 3$,

\begin{gather*}
\chi_{ed}(L(H_{1,t,t}))=
\begin{cases}
t+\ceil{\frac{t}{2}}+\ceil{\frac{t}{4}}, & t\equiv 0,2,3 \pmod{4};\\
t+\ceil{\frac{t}{2}}+\floor{\frac{t}{4}}, & t\equiv 1 \pmod{4}.\\
\end{cases}
\end{gather*}
\end{theorem}

\begin{figure}[h]
\centering
\begin{tikzpicture}[scale=0.4] 
\vertex (1) at (0*360/6:3.5) []{$c_1$};
\vertex (2) at (1*360/6:3.5) []{$c_2$};
\vertex (3) at (2*360/6:3.5) []{$c_3$};
\vertex (4) at (3*360/6:3.5) []{$c_4$};
\vertex (5) at (4*360/6:3.5) []{$c_5$};
\vertex (6) at (5*360/6:3.5) []{$c_6$};

\vertex (8) at (1*360/12:6.5) []{$c_{7}$};
\vertex (9) at (2*360/12:6.5) []{$c_{3}$};
\vertex (10) at (3*360/12:6.5) []{$c_{10}$};
\vertex (11) at (4*360/12:6.5) []{$c_{4}$};
\vertex (12) at (5*360/12:6.5) []{$c_{8}$};
\vertex (13) at (6*360/12:6.5) []{$c_{5}$};
\vertex (14) at (7*360/12:6.5) []{$c_{11}$};
\vertex (15) at (8*360/12:6.5) []{$c_{6}$};
\vertex (16) at (9*360/12:6.5) []{$c_{9}$};
\vertex (17) at (10*360/12:6.5) []{$c_{1}$};
\vertex (18) at (11*360/12:6.5) []{$c_{10}$};
\vertex (19) at (12*360/12:6.5) []{$c_{2}$};
\path
(1) edge (2)
(2) edge (3)
(3) edge (4)
(4) edge (5)
(5) edge (6)
(6) edge (1)
(1) edge (3)
(1) edge (4)
(1) edge (5)
(2) edge (4)
(2) edge (5)
(2) edge (6)
(3) edge (5)
(3) edge (6)
(4) edge (5)
(4) edge (6)
(9) edge (10)
(10) edge (11)
(11) edge (12)
(12) edge (13)
(13) edge (14)
(14) edge (15)
(15) edge (16)
(16) edge (17)
(17) edge (18)
(18) edge (19)
(19) edge (8)
(8) edge (9)
(10) edge [bend left=10]  (12)
(14) edge [bend right=10]  (12)
(14) edge [bend left=10]  (16)
(16) edge [bend left=10]  (18)
(18) edge [bend left=10]  (8)
(8) edge [bend left=10] (10)
(3) edge (11)
(4) edge (13)
(5) edge (15)
(6) edge (17)
(1) edge (19)
(2) edge (9)
(2) edge (8)
(2) edge (10)
(3) edge (10)
(3) edge (12)
(4) edge (12)
(4) edge (14)
(5) edge (14)
(5) edge (16)
(6) edge (16)
(6) edge (18)
(1) edge (18)
(1) edge (8)
;

\end{tikzpicture}
\caption{Equitable dominator coloring of $L(H_{1,6,6})$.}
\label{Helm}
\end{figure}
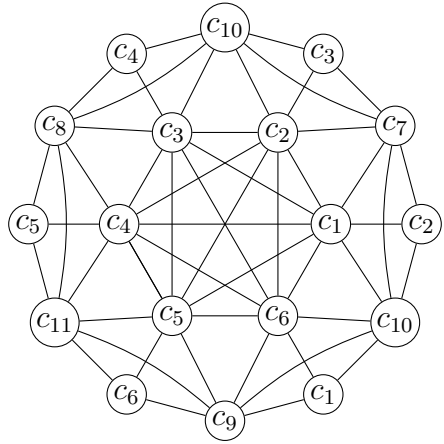

\begin{proof}
Let $V(H_{1,t,t})=\{v: deg(v)=t\}\cup\{v_i: deg(v_i)=4\} \cup \{u_i: deg(u_i)=1\}; 1 \le i \le t$ be the vertex set of $H_{1,t,t}$, and let $E(H_{1,t,t})=\{e_i: e_i=vv_i\}\cup\{e_i':e_i'=v_iv_{i+1}\}\cup\{e_i'':e_i''=v_iu_i\}$, for $1\le i\le t$ be the edge set of $H_{1,t,t}$ where the suffixes are taken under addition modulo $t$. 
Since the edges $e_i; 1\le i\le t$ have a common vertex $v$ in $H_{1,t,t}$, they induce a clique of order $t$ in $L(H_{1,t,t})$ which requires $t$ colors for its proper coloring. Based on the adjacency of the vertices $e_i; 1 \le i \le t$, the vertices $e_i'$ and $e_i''$ are assigned colors such that the vertices $e_i; 1 \le i \le t$ either dominate the color class of the color assigned to the vertices $e_i'$ or dominate the color class of the color assigned to the vertices $e_i''$, such that $1 \le i\le t$. Thus, we have the following coloring schemes.

\textit{Coloring Scheme 1:} Consider a coloring $c$ such that $c(e_i)=c_i; 1\le i\le t$ and $c(e_i'')=c(e_{i+1}); 1\le i\le t$. The remaining vertices $e_i'; 1 \le i \le t$  are assigned colors based on the value of congruence of $t$ modulo $4$ as follows.

\textit{Case 1:} In the case when $t \equiv 0,2,3 \Mod{4}$, let
\begin{gather*}
c(e_i')=
\begin{cases}
c_{t+k}, & i=2k-1, 1\le k \le \ceil{\frac{t}{2}};\\
c_{t+\ceil{\frac{t}{2}}+k}, & i\in\{i,i+2\ceil{\frac{t}{4}}\},i=2k, 1\le k\le \ceil{\frac{t}{4}}.\\
 \end{cases}
\end{gather*}

Here, the coloring $c$ of $L(H_{1,t,t})$ satisfies the criterion of dominator coloring such that when $i\equiv 1\Mod{2}$, the vertices $e_i'',e_{i+1}'',e_{i}$ dominate the color class $\{e_i'\}$, and for $1\le i\le t$, every vertex $e_i,e_i'; 1\le i \le t$ dominates the color class $V_{i+1}$. The coloring $c$ also satisfies the criterion of equitability since the cardinality of every color class $V_i$ in the coloring $c$ is at most $2$; hence $c$ is an equitable dominator coloring of $L(H_{1,t,t})$ using $t+\ceil{\frac{t}{2}}+\ceil{\frac{t}{4}}$ colors.  

\textit{Case 2:} When $t\equiv 1 \pmod{4}$, let the vertices $e_i'; 1\le i \le t$ be assigned colors such that
\begin{gather*}
c(e_i')=
\begin{cases}
c_{t+k}, & i=2k-1, 1\le k \le \ceil{\frac{t}{2}};\\
c_{t+\ceil{\frac{t}{2}}+k}, & i\in\{i,i+2\floor{\frac{t}{4}}\}, i=2k, 1\le k\le \floor{\frac{t}{4}}.\\
 \end{cases}
\end{gather*}

Following similar arguments as given in \textit{Case 1}, it can be shown that $c$ is an equitable dominator coloring of $L(H_{1,t,t})$ using $t+\ceil{\frac{t}{2}}+\floor{\frac{t}{4}}$ colors.  

\textit{Coloring Scheme 2:} Consider a coloring $c':V(L(H_{1,t,t}))\to  \mathcal{C}$, where $\mathcal{C}$ represents the set of colors used in this coloring. Let $c'(e_i)=c_{i}$ and $c'(e_i')=c(e_{i-1})$, for $1 \leq i \leq t$. Since the vertices $e_i''$ cannot dominate the color class of either the vertex $e_i,e_{i}',e_{i-1}'$: let $c'(e_{i}'')=c_{t+i}; 1\le i \le t$. Here, every vertex $e_{i}''$ dominates its own color class $\{e_i''\}$, and every vertex $e_i; 1\le i \le t$ dominates the color class $V_{i-1}$, where suffixes are taken under addition modulo $t$. Additionally, every vertex $e_i'$ dominates the color class $\{e_i,e_{i+1}'\}$. Here, the cardinality of each color class $V_i$ is at most $2$, and hence, the coloring mentioned is an equitable dominator coloring using $2t$ colors.

Consider a coloring $c''$ such that the cardinality of at least one color class is 3. In order to satisfy the equitability condition, the cardinality of all the other color classes should be at least 2 or at most 4. Due to the clique induced by the vertices $e_i; 1\le i \le t$ in $L(H_{1,t,t})$, we require at least $t$ colors for its proper coloring. Let $c''(e_1'')=c''(e_3'')=c''(e_2)=c_2$. For the vertex $e_1''$ to dominate a color class, either the vertex $e_1',e_1$ or the vertex $e_t'$ must be given a unique color such that the cardinality of the color class is 1, thus violating the equitability condition. Hence, coloring the graph $L(H_{1,t,t})$ such that the cardinality of a color class is at least 3 is not possible.

Therefore, the minimum number of colors is obtained by following the coloring as given in \textit{Coloring scheme 1}. Thus, the equitable dominator chromatic number of $L(H_{1,t,t})$ follows as given in the result.  For understanding the coloring pattern further, refer to Figure \ref{Helm}.
\end{proof}

\begin{theorem}\label{eqdomlineclosedhelm}
For $t\geq 3$, 
\begin{gather*}
\chi_{ed}(L(CH_{1,t,t})=
\begin{cases}
2t+\ceil{\frac{t}{4}}, & t\equiv 0,2\Mod{4};\\
2t+\floor{\frac{t}{4}}+1, & t\equiv 1,3\Mod{4}.\\
\end{cases}
\end{gather*}
\end{theorem}

\begin{proof}
Let $V(CH_{1,t,t})=\{v: deg(v)=t\}\cup\{v_i: deg(v_i)=4\} \cup \{u_i: deg(u_i)=3\}; 1 \le i \le t$ be the vertex set of $CH_{1,t,t}$ and let $E(CH_{1,t,t})=\{e_i: e_i=vv_i\}\cup\{e_i':e_i'=v_iv_{i+1}\}\cup\{e_i'':e_i''=v_iu_i\}\cup\{e_{i}''':e_i'''=u_iu_{i+1}\}$, for $1\le i\le t$ be the edge set of $CH_{1,t,t}$ where the suffixes are taken under addition modulo $t$. Consider a coloring $c$ of $L(CH_{1,t,t})$ such that for a vertex $w \in V(L(CH_{1,t,t}))$,
\begin{gather*}
c(w)=
\begin{cases}
c_i, & w\in\{e_i,e_{i+1}'\}, 1\le i\le t;\\
c_{t+i},& w\in\{e_i''\}, 1\le i\le t;\\
c_{2t+k}, & w\in \{e_{i}''',e_{i+\ceil{\frac{t}{2}}}'''\}, 1\le k\le \ceil{\frac{t}{2}}.\\
\end{cases}
\end{gather*}
Note that the cardinality of every color class in the coloring $c$ is at most $2$. The vertices $e_i''; 1\le i \le t$ dominate its own color class and the remaining vertices $e_i,e_{i-1}',e_i',e_{i-1}''',e_i'''$ dominate the color class $\{e_i''\}$. Thus the coloring $c$ is an equitable dominator coloring of $L(CH_{1,t,t})$ using $2t+\ceil{\frac{t}{2}}$ colors.  

Assume an alternate coloring $c'$ of $L(CH_{1,t,t})$ such that $c'(e_i)=c'(e_{i+1}'')=c_i; 1\le i \le t$. The vertices $e_i'; 1\le i \le t$ can be assigned colors such that $c'(e_i')=c'(e_{i+\ceil{\frac{t}{2}}}')=c_{t+k},1\le k\le \ceil{\frac{t}{2}}$. Here, every vertex $e_i;1\le i\le t$ dominates the color class $\{e_{i-1},e_i''\}$, and every vertex $e_i';1\le i\le t$ dominates the color class $\{e_i,e_{i+1}''\}$. The remaining vertices $e_i'''; 1\le i\le t$ are assigned colors based on the congruence of $t$ with respect to modulo $4$ as follows.

\textit{Subcase 1:} Let $t\equiv 0,2\Mod{4}$. In this case the vertices $e_i'''; 1\le i\le t$ are assigned colors as follows.
\begin{gather*}
    c'(e_i''')=
    \begin{cases}
        c_{t+\frac{t}{2}+k},& i=2k-1, 1\le k\le \frac{t}{2};\\
        c_{2t+k}, & i\in\{i,i+2\ceil{\frac{t}{4}}\}, i=2k, 1\le k\le \ceil{\frac{t}{4}}.\\
    \end{cases}
\end{gather*}

\textit{Subcase 2:} Let $t\equiv 1\Mod{4}$. In this case the vertices $e_i'''; 1\le i\le t$ are assigned colors as follows.
\begin{gather*}
    c'(e_i''')=
    \begin{cases}
        c_{t+\ceil{\frac{t}{2}}+k},& i=2k-1, 1\le k\le \ceil{\frac{t}{2}};\\
        c_{2t+1+k}, & i\in\{i,i+2\floor{\frac{t}{4}}\}, i=2k, 1\le k\le \floor{\frac{t}{4}}.\\
    \end{cases}
\end{gather*}

\textit{Subcase 3:} Let $t\equiv 3\Mod{4}$. In this case the vertices $e_i'''; 1\le i\le t$ are assigned colors as follows.
\begin{gather*}
    c'(e_i''')=
    \begin{cases}
        c_{t+\ceil{\frac{t}{2}}+k},& i=2k-1, 1\le k\le \ceil{\frac{t}{2}};\\
        c_{2t+1+k}, & i\in\{i,i+2\floor{\frac{t}{4}}\}, i=2k, 1\le k\le \floor{\frac{t}{4}};\\
        c'(e_{\ceil{\frac{t}{2}}}'),& i=t-1.\\
    \end{cases}
\end{gather*}
In the coloring $c'$, the vertices $e_i''',e_i''$ such that $i\equiv 1\Mod{2}$ dominates the color class $\{e_i'''\}$, and the vertices $e_i''$ such that $i\equiv 0\Mod{2}$ dominates the color class $\{e_{i-1}'''\}$. Thus every vertex in $V(L(CH_{1,t,t}))$ satisfies the criterion of dominator coloring in the coloring $c'$ of $L(CH_{1,t,t})$. Also note that the cardinality of every color class in the coloring  $c'$ of $L(CH_{1,t,t})$ is at most $2$. Hence $c'$ is an equitable dominator coloring of $L(CH_{1,t,t})$.

From the above mentioned colorings, the coloring $c'$ of $L(CH_{1,t,t})$ gives the minimum number of colors used in an equitable dominator coloring of $L(CH_{1,t,t})$.

Using similar arguments as given in the Proof of Theorem $2.4$, it can be observed that $L(CH_{1,t,t})$ does not admit a coloring $c''$ such that there exists a color class of cardinality greater than or equal to $3$. This proves the result.
\end{proof}

An edge $uv$ of a graph $G$ is said to be \emph{subdivided} when the edge $uv$ is deleted and a new vertex $w$ is added in between the existing vertices $u$ and $v$ of the graph $G$ such that two new edges $uw$ and $wv$ are formed. A \emph{gear graph}, denoted by $G_{1,t}$, is a graph obtained by subdividing every edge $v_{i}v_{i+1}$ of a wheel graph $W_{1,t}$.

\begin{theorem}
For $t \geq 3$,
\begin{gather*}
\chi_{ed}(L(G_{1,t}))=
\begin{cases}
\frac{7t}{4}, & t \equiv 0 \pmod{4};\\
t+\floor{\frac{t}{2}}+\floor{\frac{t}{4}}+1, & t\equiv 1,3 \pmod{4};\\
\frac{3t}{2}+\floor{\frac{t}{4}}+1, & t\equiv 2 \pmod{4}.\\
 \end{cases}
\end{gather*}
\end{theorem}

\begin{proof}
For $1\le i \le t$, let $V(G_{1,t})=\{v: deg(v)=t\}\cup \{v_i: deg(v_i)=3\} \cup \{v_i': deg(v_i')=2\}$ such that the vertex $v_i'$ is the vertex subdividing the edge $v_iv_{i+1}$ of $W_{1,t}$. Let $V(L(G_{1,t}))=\{e_i:e_i=vv_i, 1\le i\le t\} \cup \{e_i': 1 \le i \le 2t\}$, such that for $e_{2i-1}'$ is the edge between the vertices $v_i,v_i'$ and for $e_{2i}'$ is the edge between the vertices $v_i',v_{i+1}$. Note that the suffixes for $v_i,v_i',e_i$ are taken under addition modulo $t$, whereas the suffixes for $e_i'$ are taken under addition modulo $2t$. In $L(G_{1,t})$, the vertices $e_i; 1\le i \le t$ induce a clique of order $t$ which requires $t$ colors for its proper coloring. Let $c(e_i)=c_i;1\le i \le t$. Observe that in order for the vertices $e_i; 1\le i \le t$ to satisfy the criterion of dominator coloring, every vertex $e_i$ can either dominate the color class $\{e_{i+1},e_{2(i-1)}'\}$ or the color class $\{e_{i+1},e_{2(i-1)+1}'\}$. Taking this consideration into account, and  depending on the congruence of $t$ ,  the following coloring patterns can be considered.

\textit{Case 1:} Let $t \equiv 0 \pmod{4}$.

\textit{Subcase 1.1:} Consider a coloring $c$ of $L(G_{1,t})$, such that for a vertex $w \in V(L(G_{1,t}))$,
\begin{gather*}
    c(w)=
    \begin{cases}
        c_i, & w \in \{e_i,e_{2i}'\}, 1\le i\le t;\\
        c_{t+k+1},& w\in\{e_i'\},i=4k+1, \textit{ and } 0\le k \le \frac{t}{2}-1;\\
        c_{\frac{3t}{2}+k+1},&  w\in\{e_i',e_{i+t}'\},i=4k+3, \textit{ and } 0\le k \le \frac{t}{4}-1.\\
    \end{cases}
\end{gather*}

The cardinality of every color class in this coloring $c$ of $L(G_{1,t})$ is at most $2$. Here, the vertices $e_i; 1\le i \le t$ dominate the color class $V_{i-1}$, the vertices $e_i';i=2k-1 $ dominates the color class $V_k$ and the vertices $e_{i-1}',e_i',e_{i+1}'$ dominate the color class $\{e_i'\}$ such that $i\equiv 1 \Mod{4}$. Hence, the coloring $c$ is an equitable dominator coloring of $L(G_{1,t})$ using $\frac{7t}{4}$ colors. 

\textit{Subcase 1.2:} A coloring $c'$ similar to the coloring $c$ of $L(G_{1,t})$ can be considered, such that for a vertex $w \in V(L(G_{1,t}))$,
\begin{gather*}
    c'(w)=
    \begin{cases}
        c_i, & w \in \{e_i,e_{2i}'\}, 1\le i\le t;\\
        c_{t+k+1},& w\in\{e_i'\},i=4k+3,\textit{ and } 0\le k \le \frac{t}{2}-1;\\
        c_{\frac{3t}{2}+k+1},&  w\in\{e_i',e_{i+t}'\},i=4k+1,\textit{ and } 0\le k \le \frac{t}{4}-1.\\
    \end{cases}
\end{gather*}

Similar to \textit{Subcase 1.1}, the vertices $e_i; 1\le i \le t$ dominate color class $V_{i-1}$, the vertices $e_i'$ such that $i=2k+1$ dominates the color class $V_{k}$, and the vertices $e_{i-1}',e_i',e_{i+1}'$ dominates the color class $\{e_i'\}$ such that $i\equiv 3\Mod{4}$. The cardinality of every color class used in the coloring $c$ is of cardinality at most $2$. Hence, the coloring $c$ is an equitable dominator coloring using $\frac{7t}{4}$ colors.  

\textit{Subcase 1.3:} Consider a coloring $c''$ of $L(G_{1,t})$, such that for a vertex $w \in V(L(G_{1,t}))$,
\begin{gather*}
    c''(w)=
    \begin{cases}
        c_i, & w \in \{e_i,e_{2i+1}'\}, 1\le i\le t;\\
        c_{t+k}, & w\in\{e_i'\}, i=2k, \textit{ and } 1\le k\le t.\\
    \end{cases}
\end{gather*}

The coloring $c''$ satisfies the criterion of equitable coloring since the cardinality of every color class $V_i; 1\le i \le 2t$ is at most $2$. Similar to \textit{Subcase 1.1} the vertices $e_i; 1\le i\le t$ dominate the color class $V_{i-1}$. Based on the assignment of colors to the vertices $e_i,e_{2i+1}'$, it is observed that in order for the vertices $e_i'; 1\le i \le 2t$ to satisfy the criterion of dominator coloring the vertices $e_i': i\equiv 0\Mod{2}$ must be assigned unique colors such that the vertices $e_i',e_{i+1}'$, where $i\equiv 0\Mod{2},1\le i\le 2t$ dominate the color class $\{e_i'\}$. Thus, the coloring $c''$ is an equitable dominator coloring of $L(G_{1,t})$ using $2t$ colors.

\textit{Case 2:} Let $t \equiv 1 \pmod{4}$. Here, we have the following sub-cases.

\textit{Subcase 2.1:} Similar to \textit{Subcase 1.1}, consider a coloring $c$ of $L(G_{1,t})$, such that for a vertex $w \in V(L(G_{1,t}))$,
\begin{gather*}
    c(w)=
    \begin{cases}
        c_i, & w \in \{e_i,e_{2i}'\}, 1\le i\le t;\\
        c_{t+k+1},& w\in\{e_i'\},i=4k+1, \textit{ and } 0\le k \le \floor{\frac{t}{2}};\\
        c_{t+\floor{\frac{t}{2}}+k+2},&  w\in\{e_i',e_{i+4\floor{\frac{t}{4}}}'\},i=4k+3, \textit{ and } 0\le k \le \floor{\frac{t}{4}}-1.\\
    \end{cases}
\end{gather*}
Following similar arguments as given in the \textit{Subcase 1.1}, it can be proved that $c$ is an equitable dominator coloring of $L(G_{1,t})$ using $t+\floor{\frac{t}{2}}+\floor{\frac{t}{4}}+1$ colors.

\textit{Subcase 2.2:} Similar to \textit{Subcase 1.2}, consider a coloring $c'$ of $L(G_{1,t})$, such that for a vertex $w \in V(L(G_{1,t}))$,
\begin{gather*}
    c'(w)=
    \begin{cases}
        c_i, & w \in \{e_i,e_{2i}'\}, 1\le i\le t;\\
        c_{t+k+1},& w\in\{e_i'\},i=4k+3, \textit{ and } 0\le k \le \floor{\frac{t}{2}}-1;\\
        c_{t+\floor{\frac{t}{2}}+k+1},&  w\in\{e_i',e_{i+4\ceil{\frac{t}{4}}}'\},i=4k+1, \textit{ and } 0\le k \le \floor{\frac{t}{4}}.\\
    \end{cases}
\end{gather*}
Following similar arguments as given in the \textit{Subcase 1.2}, it can be proved that $c'$ is an equitable dominator coloring of $L(G_{1,t})$ using $t+\floor{\frac{t}{2}}+\floor{\frac{t}{4}}+1$ colors.

\textit{Case 3:} Let $t \equiv 2 \pmod{4}$. Here, we have the following sub-cases.

\textit{Subcase 3.1:} Consider a coloring $c$ of $L(G_{1,t})$, such that for a vertex $w \in V(L(G_{1,t}))$,
\begin{gather*}
    c(w)=
    \begin{cases}
        c_i, & w \in \{e_i,e_{2i}'\}, 1\le i\le t;\\
        c_{t+k+1},& w\in\{e_i'\},i=4k+3, \textit{ and } 0\le k \le \frac{t}{2}-1;\\
        c_{\frac{3t}{2}+k+1},&  w\in\{e_i',e_{i+4\ceil{\frac{t}{4}}}'\},i=4k+1, \textit{ and } 0\le k \le \floor{\frac{t}{4}}.\\
    \end{cases}
\end{gather*}
Following similar arguments as given in the \textit{Subcase 1.1}, it can be proved that $c$ is an equitable dominator coloring of $L(G_{1,t})$ using $\frac{3t}{2}+\floor{\frac{t}{4}}+1$ colors.

\textit{Subcase 3.2:} Consider a coloring $c'$ of $L(G_{1,t})$, such that for a vertex $w \in V(L(G_{1,t}))$,
\begin{gather*}
    c'(w)=
    \begin{cases}
        c_i, & w \in \{e_i,e_{2i}'\}, 1\le i\le t;\\
        c_{t+k+1},& w\in\{e_i'\},i=4k+1, \textit{ and } 0\le k \le \frac{t}{2}-1;\\
        c_{\frac{3t}{2}+k+1},&  w\in\{e_i',e_{i+4\ceil{\frac{t}{4}}}'\},i=4k+3, \textit{ and } 0\le k \le \floor{\frac{t}{4}}.\\
    \end{cases}
\end{gather*}
Following similar arguments as given in the \textit{Subcase 1.2}, it can be proved that $c'$ is an equitable dominator coloring of $L(G_{1,t})$  using $\frac{3t}{2}+\floor{\frac{t}{4}}+1$ colors.

\textit{Case 4:} Let $t \equiv 3 \pmod{4}$. 

\textit{Subcase 4.1:} Similar to \textit{Subcase 2.1}, consider a coloring $c$ of $L(G_{1,t})$, such that for a vertex $w \in V(L(G_{1,t}))$,
\begin{gather*}
    c(w)=
    \begin{cases}
        c_i, & w \in \{e_i,e_{2i}'\}, 1\le i\le t;\\
        c_{t+k+1},& w\in\{e_i'\},i=4k+1, \textit{ and } 0\le k \le \floor{\frac{t}{2}};\\
        c_{t+\floor{\frac{t}{2}}+k+2},&  w\in\{e_i',e_{i+4\ceil{\frac{t}{4}}}'\},i=4k+3, \textit{ and } 0\le k \le \floor{\frac{t}{4}}.\\
    \end{cases}
\end{gather*}
Following similar arguments as given in the \textit{Subcase 1.1}, it can be proved that $c$ is an equitable dominator coloring of $L(G_{1,t})$ using $t+\floor{\frac{t}{2}}+\floor{\frac{t}{4}}+2$ colors.

\textit{Subcase 4.2:} Similar to \textit{Subcase 2.2}, consider a coloring $c'$ of $L(G_{1,t})$, such that for a vertex $w \in V(L(G_{1,t}))$,
\begin{gather*}
    c'(w)=
    \begin{cases}
        c_i, & w \in \{e_i,e_{2i}'\}, 1\le i\le t;\\
        c_{t+k+1},& w\in\{e_i'\},i=4k+3, \textit{ and } 0\le k \le \floor{\frac{t}{2}}-1;\\
        c_{t+\floor{\frac{t}{2}}+k+1},&  w\in\{e_i',e_{i+4\ceil{\frac{t}{4}}}'\},i=4k+1, \textit{ and } 0\le k \le \floor{\frac{t}{4}}.\\
    \end{cases}
\end{gather*}
Following similar arguments as given in the \textit{Subcase 1.2}, it can be proved that $c'$ is an equitable dominator coloring of $L(G_{1,t})$ using $t+\floor{\frac{t}{2}}+\floor{\frac{t}{4}}+1$ colors.

Note that in the case when $t\equiv 3\Mod{4}$, the equitable dominator coloring mentioned in Subcase $4.2$ gives the minimum number of colors and hence the result.

Assume there exists a coloring $c^{*}$ such that there exists a color class of cardinality $3$ in the coloring $c^{*}$. To satisfy the criterion of equitable coloring, the cardinality of the remaining color classes should be either $2$ or $4$. Due to the clique induced by the vertices $e_i; 1\le i \le t$ of $L(G_{1,t})$, it can be observed that at least $t$ colors are required in the $\chi_{ed}$-coloring of $L(G_{1,t})$. Let the vertices $e_1,e_2',e_4'$ belong to the color class of cardinality $3$ in the coloring $c^{*}$. From the structure of the graph $L(G_{1,t})$, it can be observed that in order for the vertex $e_3'$ to satisfy the criterion of dominator coloring, the only possibility is for the vertex $e_3'$ to dominate its own color class, or for the vertex $e_3'$ to dominate the color class $\{e_2\}$. Hence, a coloring $c^{*}$ with cardinality of a color class greater than or equal to $3$ is not possible. This completes the proof.

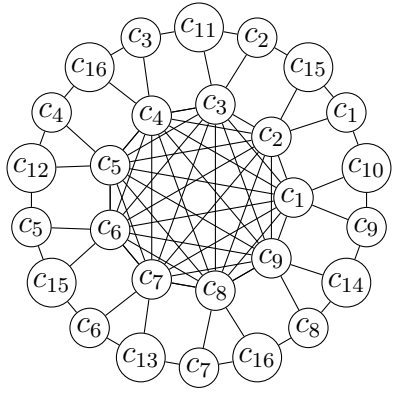
\begin{figure}[h]
\centering
\begin{tikzpicture}[scale=0.5] 
\vertex (1) at (0*360/9:2.5) []{$c_1$};
\vertex (2) at (1*360/9:2.5) []{$c_2$};
\vertex (3) at (2*360/9:2.5) []{$c_3$};
\vertex (4) at (3*360/9:2.5) []{$c_4$};
\vertex (5) at (4*360/9:2.5) []{$c_5$};
\vertex (6) at (5*360/9:2.5) []{$c_6$};
\vertex (7) at (6*360/9:2.5) []{$c_7$};
\vertex (8) at (7*360/9:2.5) []{$c_8$};
\vertex (9) at (8*360/9:2.5) []{$c_9$};
\vertex (11) at (1*360/36:4.5) []{$c_{10}$};
\vertex (12) at (3*360/36:4.5) []{$c_1$};
\vertex (13) at (5*360/36:4.5) []{$c_{15}$};
\vertex (14) at (7*360/36:4.5) []{$c_2$};
\vertex (15) at (9*360/36:4.5) []{$c_{11}$};
\vertex (16) at (11*360/36:4.5) []{$c_3$};
\vertex (17) at (13*360/36:4.5) []{$c_{16}$};
\vertex (18) at (15*360/36:4.5) []{$c_4$};
\vertex (19) at (17*360/36:4.5) []{$c_{12}$};
\vertex (20) at (19*360/36:4.5) []{$c_5$};
\vertex (21) at (21*360/36:4.5) []{$c_{15}$};
\vertex (22) at (23*360/36:4.5) []{$c_6$};
\vertex (23) at (25*360/36:4.5) []{$c_{13}$};
\vertex (24) at (27*360/36:4.5) []{$c_{7}$};
\vertex (25) at (29*360/36:4.5) []{$c_{16}$};
\vertex (26) at (31*360/36:4.5) []{$c_8$};
\vertex (27) at (33*360/36:4.5) []{$c_{14}$};
\vertex (28) at (35*360/36:4.5) []{$c_9$};
\path
(1) edge (2)
(2) edge (3)
(3) edge (4)
(4) edge (5)
(5) edge (6)
(6) edge (7)
(7) edge (8)
(8) edge (9)
(1) edge (3)
(1) edge (4)
(1) edge (5)
(1) edge (6)
(1) edge (7)
(1) edge (8)
(1) edge (9)
(2) edge (8)
(2) edge (4)
(2) edge (5)
(2) edge (6)
(2) edge (7)
(2) edge (9)
(3) edge (4)
(3) edge (5)
(3) edge (6)
(3) edge (7)
(3) edge (8)
(3) edge (9)
(4) edge (5)
(4) edge (6)
(4) edge (7)
(4) edge (8)
(4) edge (9)
(5) edge (6)
(7) edge (5)
(5) edge (8)
(5) edge (9)
(6) edge (7)
(6) edge (8)
(6) edge (9)
(7) edge (8)
(7) edge (9)
(8) edge (9)
(11) edge (28)
(1) edge (11)
(1) edge (28)
(2) edge (12)
(2) edge (13)
(3) edge (14)
(3) edge (15)
(4) edge (16)
(4) edge (17)
(5) edge (18)
(5) edge (19)
(6) edge (20)
(6) edge (21)
(7) edge (22)
(7) edge (23)
(8) edge (24)
(8) edge (25)
(9) edge (26)
(9) edge (27)
(11) edge (12)
(12) edge (13)
(13) edge (14)
(14) edge (15)
(15) edge (16)
(16) edge (17)
(17) edge (18)
(18) edge (19)
(19) edge (20)
(20) edge (21)
(21) edge (22)
(22) edge (23)
(23) edge (24)
(24) edge (25)
(25) edge (26)
(26) edge (27)
(27) edge (28);
\end{tikzpicture}
\caption{Equitable dominator coloring of $L(G_{1,9})$}
\label{gear}
\end{figure}
\end{proof}

A \emph{$t$-sunlet graph}, denoted by $Sl_t$, is a graph obtained by attaching a pendant edge to each vertex of a cycle graph $C_t$. In other words, a $t$-sunlet graph is a graph on $2t$ vertices which is obtained by taking the corona product $C_t\circ K_1$. Note that $L(Sl_t)\cong M(C_t)$ and hence we recall the following theorem.

\begin{theorem}\label{middlecycleref}
For $t\geq 3$, $\gamma(M(C_t))=\ceil{\frac{t}{2}}$.
\end{theorem}

\begin{theorem}
For $t\geq 3$, $\chi_{ed}(L(Sl_t))=t+\ceil{\frac{t}{4}}$.
\end{theorem}
\begin{proof}
Let $V(Sl_t)=\{v_i: deg(v_i)=3\} \cup \{u_i: deg(u_i)=1\}$, and let $E(Sl_t)=\{e_i: e_i=v_iv_{i+1}\} \cup \{e_i' :e_i'=v_iu_i\}$, where $1\le i\le t$, and all the suffixes are taken under addition modulo $t$. 
Note that the set $S=\{e_i: i \equiv 1 \Mod 2\}$ forms a dominating set of $L(Sl_t)$ such that $\abs{S}=\ceil{\frac{t}{2}}$ (see Theorem \ref{middlecycleref}). Here, the vertices $e_i'; 1\le i \le t$ can satisfy the criterion of dominator coloring by either dominating the color class $\{e_j\}$, such that $j=i,i+1$, or the color class $\{e_i'\}$. Thus, any equitable dominator coloring of $L(Sl_t)$ can have cardinality of every color class at most $2$. Depending on the aforementioned condition and the value of $\rho$ such that $t\equiv \rho \Mod{4}$, the following cases are considered.

\textit{Case 1:} Let $e_i'$ dominate the color class $\{e_j\}$, such that $j=i$ or $j=i+1$.

\textit{Subcase 1.1:} Let $t\equiv 0 \Mod{4}$. Consider a coloring $c$ of $L(Sl_t)$ such that for any vertex $w\in V(L(Sl_t))$,
\begin{gather*}
c(w)=
\begin{cases}
c_i,& w\in S;\\
c_{\ceil{\frac{t}{2}}+k},& w\in \{e_i',e_{i+1}'\},i=2k-1,1\le k\le \frac{t}{2};\\
c_{t+k},& w\in\{e_i,e_{i+2\floor{\frac{t}{4}}}\},i=2k,1\le k\le \frac{t}{4}.\\
\end{cases}
\end{gather*}

Here, in this coloring $c$, the vertices belonging to the dominating set of $L(Sl_t)$ are assigned unique colors such that they dominate their own color classes and all the remaining vertices in $V(L(Sl_t))-S$ dominate the color class of the vertex in $S$ belonging to their neighbourhood. The cardinality of every color class is at most $2$, and hence the coloring $c$ is an equitable dominator coloring such that $t+\frac{t}{4}$.

\textit{Subcase 1.2:} Let $t\equiv 1 \Mod{4}$. Consider a coloring $c'$ of $L(Sl_t)$ such that for any vertex $w\in V(L(Sl_t))$,
\begin{gather*}
c(w)=
\begin{cases}
c_i,& w\in S;\\
c_{\ceil{\frac{t}{2}}+k},& w\in \{e_i',e_{i+1}'\},i=2k-1,1\le k\le \ceil{\frac{t}{2}};\\
c_{t+k},& w\in\{e_i,e_{i+2\floor{\frac{t}{4}}}'\},i=2k,1\le k\le \floor{\frac{t}{4}}.\\
\end{cases}
\end{gather*}

Using similar arguments as in Subcase $1.1$, it can be established that $c'$ is an equitable dominator coloring of $L(Sl_t)$ using $t+\ceil{\frac{t}{4}}$. 

\textit{Subcase 1.3:} Let $t\equiv 2,3 \Mod{4}$. Consider a coloring $c''$ of $L(Sl_t)$ such that for any vertex $w\in V(L(Sl_t))$,
\begin{gather*}
c''(w)=
\begin{cases}
c_i,& w\in S;\\
c_{\ceil{\frac{t}{2}}+k},& w\in\{e_i,e_{i+2\ceil{\frac{t}{4}}}\},i=2k, 1\le i\le \ceil{\frac{t}{4}};\\
c_{\ceil{\frac{t}{2}}+\ceil{\frac{t}{4}}+k}, &w\in\{e_i',e_{i+1}'\},i=2k-1,1\le k\le \floor{\frac{t}{2}};\\
c(e_{2\ceil{\frac{t}{4}}}),& w=e_t'.\\
\end{cases}
\end{gather*}

Using similar arguments as in Subcase $1.1$, it can be established that $c'$ is an equitable dominator coloring of $L(Sl_t)$ using $t+\ceil{\frac{t}{4}}$. 

\textit{Case 2:} Let the vertices $e_i'; 1\le i\le t$ dominate its own color class. Consider a coloring $c'''$ of $L(Sl_t)$ such that $c(e_i')=c_i; 1\le i\le t$ and $c(e_i)=c(e_{i+\ceil{\frac{t}{2}}})=c_k$, where $1\le k\le \ceil{\frac{t}{2}}$. The coloring $c'''$ is an equitable dominator coloring of $L(Sl_t)$ using $t+\ceil{\frac{t}{2}}$ colors such that every vertex $e_i,e_i'$ dominates the color class $\{e_i'\}$, and the cardinality of every color class is at most $2$. 

From the above-mentioned colorings, it can be observed that Case $1$ gives an equitable dominator coloring of $L(Sl_t)$ using the least possible number of colors. Hence, the result follows.

\end{proof}

\begin{figure}[h]
\centering
\begin{tikzpicture}[scale=0.6] 
\vertex (1) at (0*360/5:2.5) []{$c_{1}$};
\vertex (2) at (1*360/5:2.5) []{$c_4$};
\vertex (3) at (2*360/5:2.5) []{$c_2$};
\vertex (4) at (3*360/5:2.5) []{$c_4$};
\vertex (5) at (4*360/5:2.5) []{$c_{3}$};
\vertex (6) at (1*360/10:4.5) []{$c_{5}$};
\vertex (7) at (3*360/10:4.5) []{$c_{5}$};
\vertex (8) at (5*360/10:4.5) []{$c_{6}$};
\vertex (9) at (7*360/10:4.5) []{$c_{6}$};
\vertex (10) at (9*360/10:4.5)[]{$c_{7}$};
\path
(1) edge (2)
(2) edge (3)
(3) edge (4)
(4) edge (5)
(5) edge (1)
(1) edge (6)
(6) edge (2)
(2) edge (7)
(7) edge (3)
(3) edge (8)
(8) edge (4)
(4) edge (9)
(9) edge (5)
(5) edge (10)
(10) edge (1);
\end{tikzpicture}
\caption{Equitable dominator coloring of $L(Sl_{5})$}
\label{sunlet}
\end{figure}
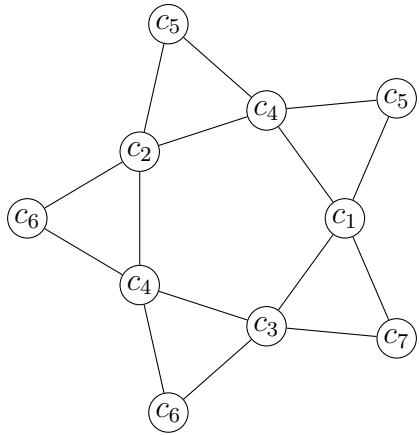

A \textit{friendship graph} denoted by $F_t$ is a graph obtained by joining $t$ copies of $C_3$ to a common vertex such that the vertex becomes the graph's universal vertex.

\begin{theorem}
 For $t \geq 2$, $\chi_{ed}(L(F_t))=2t$.    
\end{theorem}

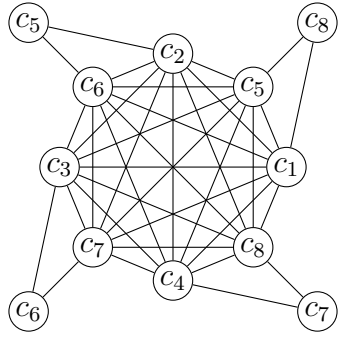
\begin{figure}[h]
\centering
\begin{tikzpicture}[scale=0.6] 
\vertex (1) at (0*360/8:2.5) []{$c_1$};
\vertex (2) at (1*360/8:2.5) []{$c_5$};
\vertex (3) at (2*360/8:2.5) []{$c_2$};
\vertex (4) at (3*360/8:2.5) []{$c_6$};
\vertex (5) at (4*360/8:2.5) []{$c_3$};
\vertex (6) at (5*360/8:2.5) []{$c_7$};
\vertex (7) at (6*360/8:2.5) []{$c_4$};
\vertex (8) at (7*360/8:2.5) []{$c_8$};
\vertex (11) at (1*360/8:4.5) []{$c_{8}$};
\vertex (12) at (3*360/8:4.5) []{$c_5$};
\vertex (13) at (5*360/8:4.5) []{$c_{6}$};
\vertex (14) at (7*360/8:4.5) []{$c_7$};
\path
(1) edge (2)
(1) edge (3)
(1) edge (4)
(1) edge (5)
(1) edge (6)
(1) edge (7)
(1) edge (8)
(2) edge (3)
(2) edge (4)
(2) edge (5)
(2) edge (6)
(2) edge (7)
(2) edge (8)
(3) edge (4)
(3) edge (5)
(3) edge (6)
(3) edge (7)
(3) edge (8)
(4) edge (5)
(4) edge (6)
(4) edge (7)
(4) edge (8)
(5) edge (6)
(5) edge (7)
(5) edge (8)
(6) edge (7)
(6) edge (8)
(7) edge (8)
(1) edge (11)
(11) edge (2)
(3) edge (12)
(12) edge (4)
(5) edge (13)
(13) edge (6)
(7) edge (14)
(14) edge (8);
\end{tikzpicture}
\caption{Equitable dominator coloring of $L(F_{4})$}
\label{friendship}
\end{figure}

\begin{proof}
For $1 \leq i \leq t$, let $u_i$ and $v_i$ be the vertices of degree $2$ such that $u_i$ and $v_i$ are the vertices belonging to the $i$-th copy of $C_3$ in $F_t$. Let the universal vertex of $F_t$ be denoted by $v$ such that $\deg(v)=2t$. Let $V(L(F_t))=\{e_i: e_i=vu_i\} \cup \{e_i': e_i'=u_iv_i\} \cup \{e_i'': e_i''=vv_i\}$, where $1 \le i \le t$, and all the suffixes are taken under addition modulo $t$. Given the edges $\{e_i:1 \le i \le t\}\cup\{e_i'': 1 \le i \le t\}$ of $F_t$ share a common vertex $v$, the vertices $\{e_i:1 \le i \le t\}\cup\{e_i'': 1 \le i \le t\}$ of $L(F_t)$ induce a clique of order $2t$. Thus, $\chi_{ed}(L(F_t))\geq 2t$. Consider a coloring $c:V(L(F_{t}))\to \{c_1,c_2,\ldots\}$, such that for a vertex $w \in V(L(F_t))$,
\begin{gather*}
c(w)=
\begin{cases}
c_i, & w\in \{e_i\}, 1\le i \le t;\\
c_{t+i}, & w\in \{e_i'',e_{i+1}'\}, 1\le i \le t.\\
\end{cases}
\end{gather*}
Refer to Figure \ref{friendship} for the coloring pattern. The coloring $c$ is an equitable dominator coloring of $L(F_t)$ since the vertices $e_i,e_i',e_i''; 1\le i \le t$ dominate the color class $V_i$, and the cardinality of every color class $V_i; 1\le i \le 2t$ in the coloring $c$ is at most $2$. Thus, $\chi_{ed}(L(F_t))\le 2t$, and the result follows.      
\end{proof}

A \textit{flower graph}, denoted by $F_{1,t}$, is a graph obtained from a helm graph $H_{1,t,t}$ by joining each of its pendant vertices to its central vertex $v$. 

\begin{theorem}
For $t\geq 3$, $\chi_{ed}(L(F_{1,t}))=2t$.    
\end{theorem}
\begin{proof}
Let $V(F_{1,t})=\{v: deg(v)=2t\} \cup \{v_i: deg(v_i)=4\} \cup \{u_i: deg(u_i)=2\}$, and let $E(F_{1,t})= \{e_i: e_i=vv_i\} \cup \{e_i': e_i'=v_iv_{i+1}\} \cup \{e_i'': e_i''=v_iu_i\} \cup \{e_i''': e_i'''=vu_i\}$, such that $1 \le i \le t$. The edges $e_i; 1 \le i \le t$ and the edges $e_i'''; 1 \le i\le t$ being adjacent to a common vertex $v$ in $F_{1,t}$ induce a clique of order $2t$ in $L(F_{1,t})$. Thus, $\chi_{ed}(L(F_{1,t}))\geq 2t$. 
Consider a coloring $c:V(L(F_{1,t}))\to \{c_1,c_2\ldots\}$, such that for $w \in V(F_{1,t})$, 
\begin{gather*}
    c(w)=
\begin{cases}
    c_i,& w\in \{e_i,e_{i+1}''\}, 1\le i \le t;\\[15pt]
    c_{t+i}, & w\in \{e_i',e_{i}'''\}, 1\le i\le t.\\
\end{cases}
\end{gather*}

Here, the suffixes are taken under addition modulo $t$. Since the cardinality of every color class $V_i; 1\le i \le 2t$ in the coloring $c$ is $2$, the coloring $c$ of $L(F_{1,t})$ satisfies the criterion of equitable coloring. The coloring $c$ of $L(F_{1,t})$ also satisfies the criterion of dominator coloring such that the vertices $e_i,e_i''; 1\le i \le t$ dominate the color class $\{e_i',e_{i}'''\}$, the vertices $e_i'; 1\le i \le t$ dominate the color class $\{e_i,e_{i+1}''\}$, and the vertices $e_i'''; 1\le i \le t$ dominate the color class $\{e_i'',e_{i-1}\}$. Hence, the coloring $c$ is an equitable dominator coloring of $L(F_{1,t})$ using $2t$ colors, and this establishes the result.
\end{proof}

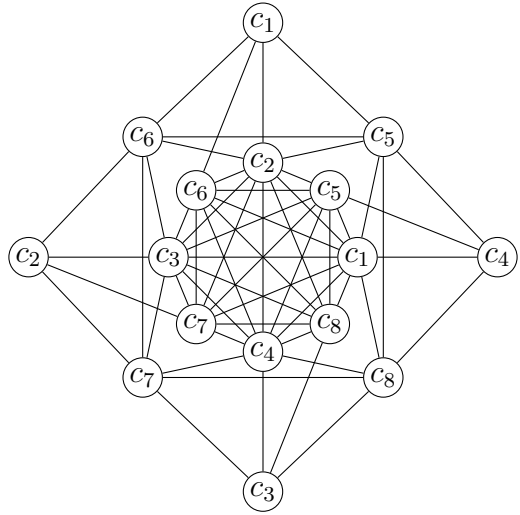
\begin{figure}[h]
\centering
\begin{tikzpicture}[scale=0.5] 
\vertex (1) at (0*360/8:2.5) []{$c_1$};
\vertex (2) at (1*360/8:2.5) []{$c_5$};
\vertex (3) at (2*360/8:2.5) []{$c_2$};
\vertex (4) at (3*360/8:2.5) []{$c_6$};
\vertex (5) at (4*360/8:2.5) []{$c_3$};
\vertex (6) at (5*360/8:2.5) []{$c_7$};
\vertex (7) at (6*360/8:2.5) []{$c_4$};
\vertex (8) at (7*360/8:2.5) []{$c_8$};
\vertex (11) at (1*360/8:4.5) []{$c_{5}$};
\vertex (12) at (3*360/8:4.5) []{$c_{6}$};
\vertex (13) at (5*360/8:4.5) []{$c_{7}$};
\vertex (14) at (7*360/8:4.5) []{$c_{8}$};
\vertex (15) at (0*360/8:6.2) []{$c_{4}$};
\vertex (16) at (2*360/8:6.2) []{$c_{1}$};
\vertex (17) at (4*360/8:6.2) []{$c_{2}$};
\vertex (18) at (6*360/8:6.2) []{$c_{3}$};

\path

(11) edge (1)
(11)edge (3)
(12) edge (3)
(12) edge (5)
(13) edge (5)
(13) edge (7)
(14) edge (7)
(14) edge (1)
(11) edge (12)
(12) edge (13)
(13) edge (14)
(14) edge (11)
(15) edge (11)
(15) edge (14)
(15) edge (1)
(15) edge (2)
(16) edge (11)
(16) edge (12)
(16) edge (4)
(16) edge (3)
(17) edge (12)
(17) edge (13)
(17) edge (5)
(17) edge (6)
(18) edge (13)
(18) edge (14)
(18) edge (7)
(18) edge (8)

(1) edge (2)
(1) edge (3)
(1) edge (4)
(1) edge (5)
(1) edge (6)
(1) edge (7)
(1) edge (8)
(2) edge (3)
(2) edge (4)
(2) edge (5)
(2) edge (6)
(2) edge (7)
(2) edge (8)
(3) edge (4)
(3) edge (5)
(3) edge (6)
(3) edge (7)
(3) edge (8)
(4) edge (5)
(4) edge (6)
(4) edge (7)
(4) edge (8)
(5) edge (6)
(5) edge (7)
(5) edge (8)
(6) edge (7)
(6) edge (8)
(7) edge (8)

;
\end{tikzpicture}
\caption{Equitable dominator coloring of $L(F_{1,4})$}
\label{flower}
\end{figure}

\textit{Double wheel} graphs denoted by $DW_{1,t}$ are graphs formed by the join of $K_1+2 C_t$.

\begin{theorem}
For $t\geq 3$, $\chi_{ed}(L(DW_{1,t}))=2t$.
\end{theorem}

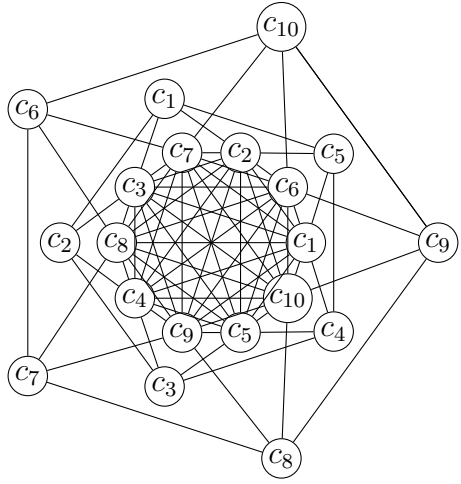
\begin{figure}[h]
\centering
\begin{tikzpicture}[scale=0.5] 
\vertex (1) at (0*360/10:2.5) []{$c_1$};
\vertex (2) at (1*360/10:2.5) []{$c_6$};
\vertex (3) at (2*360/10:2.5) []{$c_2$};
\vertex (4) at (3*360/10:2.5) []{$c_7$};
\vertex (5) at (4*360/10:2.5) []{$c_3$};
\vertex (6) at (5*360/10:2.5) []{$c_8$};
\vertex (7) at (6*360/10:2.5) []{$c_4$};
\vertex (8) at (7*360/10:2.5) []{$c_9$};
\vertex (9) at (8*360/10:2.5) []{$c_5$};
\vertex (10) at (9*360/10:2.5) []{$c_{10}$};

\vertex (11) at (1*360/10:4.0) []{$c_{5}$};
\vertex (12) at (2*360/10:6.0) []{$c_{10}$};
\vertex (13) at (3*360/10:4.0) []{$c_{1}$};
\vertex (14) at (4*360/10:6.0) []{$c_6$};
\vertex (15) at (5*360/10:4.0) []{$c_{2}$};
\vertex (16) at (6*360/10:6.0) []{$c_7$};
\vertex (17) at (7*360/10:4.0) []{$c_{3}$};
\vertex (18) at (8*360/10:6.0) []{$c_{8}$};
\vertex (19) at (9*360/10:4.0) []{$c_{4}$};
\vertex (20) at (0*360/10:6.0) []{$c_9$};

\path
(20) edge (12)
(1) edge (2)
(1) edge (3)
(1) edge (4)
(1) edge (5)
(1) edge (6)
(1) edge (7)
(1) edge (8)
(1) edge (9)
(1) edge (10)
(2) edge (3)
(2) edge (4)
(2) edge (5)
(2) edge (6)
(2) edge (7)
(2) edge (8)
(2) edge (9)
(2) edge (10)
(3) edge (4)
(3) edge (5)
(3) edge (6)
(3) edge (7)
(3) edge (8)
(3) edge (9)
(3) edge (10)
(4) edge (5)
(4) edge (6)
(4) edge (7)
(4) edge (8)
(4) edge (9)
(4) edge (10)
(5) edge (6)
(5) edge (7)
(5) edge (8)
(5) edge (9)
(5) edge (10)
(6) edge (7)
(6) edge (8)
(6) edge (9)
(6) edge (10)
(7) edge (8)
(7) edge (9)
(7) edge (10)
(8) edge (9)
(8) edge (10)
(9) edge (10)
(11) edge (1)
(11) edge (3)
(12) edge (2)
(12) edge (4)
(13) edge (3)
(13) edge (5)
(14) edge (4)
(14) edge (6)
(15) edge (5)
(15) edge (7)
(16) edge (6)
(16) edge (8)
(17) edge (7)
(17) edge (9)
(18) edge (8)
(18) edge (10)
(19) edge (9)
(19) edge (1)
(20) edge (10)
(20) edge (2)
(11) edge (13)
(13) edge (15)
(15) edge (17)
(17) edge (19)
(19) edge (11)
(12) edge (14)
(14) edge (16) 
(16) edge (18)
(18) edge (20)         
(20) edge (12)
;
\end{tikzpicture}
\caption{Equitable dominator coloring of $L(DW_{1,5})$}
\label{doublewheel}
\end{figure}

\begin{proof}
Let $V(DW_{1,t})=\{v\}\cup\{v_i: 1\le i\le t\}\cup\{u_i: 1\le i\le t\}$ such that the vertices $v_i$ form the rim vertices of the inner cycle $C_t$, and the vertices $u_i$ form the rim vertices of the outer cycle $C_t$. For $1\le i\le t$, let $V(L(DW_{1,t}))=\{e_i: e_i=vv_i\} \cup \{e_i':e_i'=v_iv_{i+1}\} \cup \{e_i'':e_i''=vu_i\}\cup \{e_i''': e_i'''=u_iu_{i+1}\}$, such that the suffixes are taken under addition modulo $t$. The edges $\{e_i: 1\le i\le t\}\cup\{e_i'': 1\le i\le t\}$ of $DW_{1,t}$ being adjacent to the vertex $v$ induce a clique of order $2t$ in $L(DW_{1,t})$, which requires $2t$ colors for its proper coloring. Thus, $\chi_{ed}(L(DW_{1,t})\geq 2t$. Consider a coloring $c:V(L(DW_{1,t}))\to \{c_1,c_2\ldots\}$, such that for $w \in V(DW_{1,t})$, 
\begin{gather*}
    c(w)=
\begin{cases}
    c_i,& w\in \{e_i,e_{i+1}'\}, 1\le i \le t;\\[15pt]
    c_{t+i}, & w\in \{e_{i}'',e_{i+1}'''\}, 1\le i\le t.\\
\end{cases}
\end{gather*}

Refer to Figure \ref{doublewheel} for the coloring pattern.  Here, the cardinality of every color class $V_i; 1\le i \le 2t$ in the coloring $c$ is $2$. Every vertex $e_i; 1\le i \le t$ of $L(DW_{1,t})$ dominates the color class $\{e_{i-1},e_i'\}$ and every vertex $e_i'; 1\le i\le t$ of $L(DW_{1,t})$ dominates the color class $\{e_i,e_{i+1}'\}$. Similarly, every vertex $e_i''; 1\le i \le t$ of $L(DW_{1,t})$ dominates the color class $\{e_{i-1}'',e_i'''\}$ and every vertex $e_i'''; 1\le i\le t$ of $L(DW_{1,t})$ dominates the color class $\{e_i'',e_{i+1}'''\}$.  Thus, the coloring mentioned is an equitable dominator coloring of $L(DW_{1,t})$, and the result follows.    
\end{proof}

\section{Conclusion}
In this paper, the concept of equitable dominator coloring is explored for line graphs of some standard graph classes such as wheel graphs, complete bipartite graphs, flower graphs, and gear graphs, and the corresponding parameter is obtained. Some open problems related to this topic are suggested below.

\begin{enumerate}[label={(\roman*)}]
    \item The equitable dominator chromatic number for the closed helm graph and its derived families can be obtained.
    \item The concept of equitable dominator coloring can be explored for various graph operations, such as vertex addition, vertex deletion, edge subdivision, and deletion of an edge.
    \item The parameter can be studied for the central graphs of various graph families.
\end{enumerate}

\noi All these highlight the broad scope of the topic for further investigation.

\nocite{*}
\bibliographystyle{abbrv}
\bibliography{ref}
\end{document}